\documentclass[11pt]{article}

\usepackage{amsmath}
\usepackage{amssymb}

\def\co{\colon\thinspace}
\newcommand{\Z}{{\mathbb Z}}
\newcommand{\Q}{{\mathbb Q}}
\newcommand{\R}{{\mathbb R}}
\newcommand{\C}{{\mathbb C}}
\newcommand{\LL}{{\mathbb L}}

\usepackage{amsthm}

\newtheorem{thm}{Theorem}

\newtheorem{prop}[thm]{Proposition}
\newtheorem{cor}[thm]{Corollary}

\theoremstyle{definition}

\newtheorem*{rems}{Remarks}

\newtheorem*{defn}{Definition}
\newtheorem*{ex}{Example}

\begin{document}

\title{Contact Dehn surgery,
symplectic fillings,\\and Property P for knots}

\author{Hansj\"org Geiges\\ 
        Mathematisches Institut, Universit\"at zu K\"oln\\
         Weyertal 86--90, 50931 K\"oln, Germany}

\date{Mathematische Arbeitstagung June 10--16, 2005\\
      Max-Planck-Institut f\"ur Mathematik, Bonn, Germany}
\maketitle

\section{Property P for knots}
According to a fundamental theorem of Lickorish
and Wallace from the 1960s, every closed, connected, orientable
$3$--manifold can be obtained by performing Dehn surgery on a link
in the $3$--sphere. Previous to the recent work of Perelman, which
is expected to close the coffin on the Poincar\'e conjecture,
it was a natural question for geometric topologists whether
one might be able to produce a counterexample to that
conjecture by a single Dehn surgery. This led to the
definition of the following property, whose name is generally regarded
as a little unfortunate.

\begin{defn}
A knot $K$ in $S^3$ has {\bf Property P} if every nontrivial surgery
along $K$ yields a non-simply-connected $3$--manifold.
\end{defn}

Our knots are always understood to be smooth, or at least tame, i.e.\
equivalent to a smooth one.

Let me briefly recall the notion of
Dehn surgery along a knot $K$ in the $3$--sphere
$S^3$. Write $\nu K\cong S^1\times D^2$
for a (closed) tubular neighbourhood of~$K$. On the boundary
$\partial (\nu K)\cong T^2$ of this tubular neighbourhood
there are two distinguished curves (which we implicitly identify with
the classes they represent in the homology group $H_1(T^2)$):

\begin{itemize}
\item[1.] The meridian $\mu$, defined as a simple closed curve that
generates the kernel of the homomorphism on $H_1$ induced by
the inclusion $T^2\rightarrow \nu K$.
\item[2.] The preferred longitude $\lambda$, defined as a simple closed curve
that generates the kernel of the homomorphism on $H_1$ induced by
the inclusion $T^2\rightarrow C:=\overline{S^3\setminus\nu K}$.
\end{itemize}

This preferred longitude can also be characterised by the property that
it has linking number zero with~$K$. The knot $K$ bounds an embedded
surface in $S^3$ (called a {\bf Seifert surface} for~$K$), and
$\lambda$ can be obtained by pushing $K$ along that surface. For that reason,
the trivialisation of the normal bundle of $K$ defined by $\lambda$
is called the {\bf surface framing} of~$K$.

Given an orientation of $S^3$, orientations of $\mu$ and $\lambda$ are
chosen such that $(\mu ,\lambda )$ is a positive basis for $H_1(T^2)$,
with $T^2$ oriented as the boundary of~$\nu K$. In the contact geometric
setting below, the orientation of $S^3$ will be the one induced from
the contact structure.

Let $p,q$ be coprime integers. The manifold $K_{p/q}$ obtained from
$S^3$ by {\bf Dehn surgery} along $K$ with {\bf surgery coefficient}
$p/q\in\Q\cup\{\infty\}$ is defined as
\[ K_{p/q}:=\overline{S^3\setminus\nu K}\cup_g S^1\times D^2,\]
where the gluing map $g$ sends the meridian $*\times\partial D^2$
to $p\mu+q\lambda$. The resulting manifold is completely determined
by the knot and the surgery coefficient.

A simple Mayer-Vietoris argument shows that $H_1(K_{p/q})\cong\Z_{|p|}$.
Therefore, saying that a knot $K$ has Property P is equivalent to
\[\pi_1(K_{1/q})=1\;\;\;\mbox{\rm only for}\;\;\; q=0.\]
(Observe that $p/q=\infty$ corresponds to a trivial surgery.)

\begin{ex}
The unknot does {\em not} have Property P. Indeed,
every $(1/q)$--surgery on the unknot yields~$S^3$, which is seen as follows.
If $K$ is the unknot, then the closure $C$ of $S^3\setminus \nu K$
is also a solid torus. Write $\mu_C$ and $\lambda_C$ for meridian
and preferred longitude on~$\partial C$. We may
assume $\mu =\lambda_C$ and $\lambda =\mu_C$. When performing
$(1/q)$--surgery on $K$, a solid torus is glued to $C$ by sending
its meridian $\mu_0$ to $\mu +q\lambda =\lambda_C+q\mu_C$. Now, there clearly
is a diffeomorphism of $C$ that sends $\mu_C$ to itself and $\lambda_C$
to $\lambda_C+q\mu_C$. It follows that the described surgery is equivalent
to the one where we send $\mu_0$ to $\lambda_C=\mu$, which is a trivial
$\infty$--surgery.
\end{ex}

In the early 1970s,  Bing and Martin, as well as
Gonz\'alez-Acu\~na, conjectured that every
nontrivial knot has Property~P. By work of Kronheimer and
Mrowka~\cite{krmr04}, this is now a theorem.

\begin{thm}[Kronheimer-Mrowka]
\label{thm:PP}
Every nontrivial knot in $S^3$ has Property~P.
\end{thm}

Before describing the role that contact geometry has played in the
proof of this theorem, I want to indicate the importance of this
theorem beyond the negative statement that counterexamples to the
Poincar\'e conjecture cannot result from a single surgery.

\begin{prop}
If two knots $K,K'$ in $S^3$ have homeomorphic complements and one of
the knots has property~P, then the knots are equivalent, i.e.\
there is a homeo\-morphism of $S^3$ mapping $K$ to~$K'$.
\end{prop}

\begin{proof}
According to a result of Edwards~\cite{edwa64}, two compact
$3$--manifolds with boundary are homeomorphic if and only if
their interiors are homeomorphic. Thus, if $S^3\setminus K$
is homeomorphic to $S^3\setminus K'$, then there is
a homeomorphism $\varphi\co C\rightarrow C'$,
where $C:=\overline{S^3\setminus\nu K}$
and $C':=\overline{S^3\setminus\nu K'}$.

Suppose $K$ has Property~P. This implies that there is a unique
way of attaching a solid torus $S^1\times D^2$ to $C$ such that
the resulting manifold is the $3$--sphere. Hence $\varphi$ extends
to a homeomorphism $S^3\rightarrow S^3$, i.e.\ the knots $K$ and $K'$
are equivalent.
\end{proof}

Observe that in this proof we only used the weaker property that
nontrivial surgery along $K$ does not yield the standard $3$--sphere.
This had been proved earlier (for $K$ different from the unknot) by
Gordon and Luecke~\cite{golu89}. Since the unknot is characterised
by its complement being a solid torus, the result of Kronheimer and
Mrowka (or the weaker one by Gordon and Luecke) yields the following
corollary.

\begin{cor}
If two knots in $S^3$ have homeomorphic
complements, then the knots are equivalent.\hfill\qed
\end{cor}

Of course, together with a positive answer to the Poincar\'e conjecture,
the result of Gordon-Luecke implies that of Kronheimer-Mrowka.
\section{Contact Dehn surgery}
This section gives a brief report on joint work with
Fan Ding~\cite{dige04}. Recall that a (coorientable) {\bf contact structure}
$\xi$ on a differential $3$--manifold is a tangent $2$--plane field defined
as the kernel of a global differential $1$--form $\alpha$ that satisfies
the nonintegrability condition $\alpha\wedge d\alpha\neq 0$
(meaning that $\alpha\wedge d\alpha$
vanishes nowhere). An example is the standard contact
structure
\[ \xi_{st}=\ker (x\, dy-y\, dx+z\, dt-t\, dz)\]
on $S^3\subset\R^4$. This can also be characterised as the complex line
in the tangent bundle $TS^3$ of $S^3$ with respect to complex multiplication
induced from the inclusion $TS^3\subset T\C^2|_{S^3}$.

I shall have to use a few notions from contact geometry without
time for much explanation (tight and overtwisted contact structures,
convex surfaces in contact $3$--manifolds). For more details
see the introductory lectures by Etnyre~\cite{etny03} or the
{\em Handbook} chapter by the present author~\cite{geig05}.

A (smooth) knot $K$ in a contact $3$--manifold $(M,\xi )$ is called
{\bf Legendrian} if it is everywhere tangent to~$\xi$. The normal
bundle of such a knot has a canonical trivialisation, determined by
a vector field along $K$ that is everywhere transverse to~$\xi$. This
will be referred to as the {\bf contact framing}. We now consider
Dehn surgery along $K$ with coefficient $p/q$ as before, but we define
the surgery coefficient with respect to the contact framing.

It turns out that for $p\neq 0$ one can always extend the contact
structure $\xi|_{M\setminus\nu K}$ to one on the surgered manifold
in such a way
that the extended contact structure is tight on the glued-in solid torus
$S^1\times D^2$. Moreover, subject to this tightness condition there
are but finitely many choices for such an extension, and for $p/q=1/k$
with $k\in\Z$ the extension is in fact unique. These observations hinge on
the fact that $\partial (\nu K)$ is a convex surface, i.e.\ a surface admitting
a transverse flow preserving the contact structure. On solid tori with
convex boundary condition, tight contact structures have been classified
by Giroux and Honda. Furthermore, one knows how to glue contact
manifolds along convex surfaces, since the germ of a contact structure
along a convex surface is determined by some simple data on that surface.

We can therefore speak sensibly of {\bf contact $(1/k)$--surgery}.
The following theorem is proved in~\cite{dige04}.

\begin{thm}
\label{thm:dige}
Let $(M,\xi )$ be a closed, connected contact $3$--manifold. Then
$(M,\xi )$ can be obtained from $(S^3,\xi_{st})$ by contact $(\pm 1)$--surgery
along a Legendrian link.\hfill\qed
\end{thm}

\begin{rems}
(1) There is a related theorem, due to Lutz-Martinet in the early
1970s, cf.~\cite{geig05}, saying that every (closed, orientable)
$3$--manifold admits a contact structure in each homotopy class of
tangent $2$--plane fields. The original proof is based on surgery
along a link in $S^3$ {\em transverse} to $\xi_{st}$. For an alternative
proof using Legendrian surgery see~\cite{dgs04}.

(2) From the topological point of view, surgeries with integer surgery
coefficient are best, since they correspond to attaching $2$--handles
to the boundary of a $4$--manifold. Thus, contact $(\pm 1)$--surgeries
are best from both the topological and contact geometric viewpoint.

(3) If $(M',\xi')$ is obtained from $(M,\xi )$ by a contact
$(1/k)$--surgery, one can recover $(M,\xi )$ by a suitable contact
$(-1/k)$--surgery on $(M',\xi')$, see~\cite{dige04}.

(4) Contact $(-1)$--surgery is symplectic handlebody surgery in the sense
of Eliashberg and Weinstein, cf.~\cite{dgs04},
and preserves the property of being
strongly symplectically fillable (see below).
\end{rems}
\section{Symplectic fillings}
Contact geometry enters the proof of Theorem~\ref{thm:PP}
via the notion of symplectic fillings.
Observe that a contact $3$--manifold $(M,\xi )$
is naturally oriented --- the sign of the volume form
$\alpha\wedge d\alpha$ does not depend on the choice of $1$--form $\alpha$
defining a given~$\xi$; similarly, a symplectic $4$--manifold $(W,\omega )$,
i.e.\ with $\omega$ a closed $2$--form satisfying $\omega^2\neq 0$,
is naturally oriented by the volume form~$\omega^2$.

\begin{defn}
\label{defn:filling}
(a) A compact symplectic $4$--manifold $(W,\omega )$
is called a {\bf weak (symplectic) filling} of the contact manifold
$(M,\xi )$ if $\partial W=M$ as oriented manifolds (outward normal
followed by orientation of $M$ gives orientation of~$W$) and
$\omega|_{\xi}\neq 0$.

(b) A compact symplectic $4$--manifold $(W,\omega )$
is called a {\bf strong (symplectic) filling} of the contact manifold
$(M,\xi )$ if $\partial W=M$ and there is a Liouville vector
field $X$ defined near $\partial W$, pointing outwards
along $\partial W$, and satisfying $\xi =\ker (i_X\omega |_{TM})$.
Here {\bf Liouville vector field} means that the Lie derivative
${\mathcal L}_X\omega$, which is the same as $d(i_X\omega )$ because of
$d\omega =0$ and Cartan's formula, is required to equal~$\omega$.
\end{defn}

For instance, $(S^3,\xi_{st})$ is strongly filled by the standard symplectic
$4$--disc $D^4$ with $\omega_{st} =dx\wedge dy+dz\wedge dt$.
The Liouville vector
field here is the radial vector field $X=r\partial_r/2$.

It is clear that every strong filling is also a weak filling. The converse
is false: There are contact structures that are weakly but not strongly
fillable; such examples are due to Eliashberg and Ding-Geiges.

The contact geometric result that allowed Kronheimer and Mrowka to conclude
their proof of Property P was established by Eliashberg~\cite{elia04}.

\begin{thm}[Eliashberg]
\label{thm:elia-filling}
Any weak symplectic filling of a contact $3$--manifold embeds
symplectically into a closed symplectic $4$--manifold.
\end{thm}

An alternative proof was given by Etnyre~\cite{etny04}. Both proofs rely
on open book decompositions adapted (in the sense of
Giroux) to contact structures.
Theorem~\ref{thm:elia-filling} being a cobordism
theoretic result, it is arguably more natural to give a surgical proof.
\"Ozba\u{g}c\i\ and Stipsicz~\cite{ozst04} were the first to observe that
such a proof, based on Theorem~\ref{thm:dige}, can indeed be devised.
In the remainder of this section, I shall sketch this surgical argument.

\vspace{2mm}

Theorem~\ref{thm:elia-filling} is proved by showing that any contact
$3$--manifold has what is
called a {\bf concave} filling that can be glued to the given
(convex) filling. (For instance, a strong concave filling
corresponds to a Liouville vector field pointing inwards
along the boundary.) Such a `cap', attached to the
(convex) symplectic filling of the contact manifold, gives the desired
closed symplectic manifold.

\vspace{2mm}

(i) Strong fillings can be capped off: Let $(W,\omega )$ be a strong
filling of $(M,\xi )$. By Theorem~\ref{thm:dige}, there
is a Legendrian link $\LL =\LL^-\sqcup\LL^+$ in $(S^3,\xi_{st})$
such that contact $(-1)$--surgery along the components of $\LL^-$
and contact $(+1)$--surgery along those of $\LL^+$ produces $(M,\xi )$.
By Remarks (3) and (4) we can attach symplectic $2$--handles to
the boundary $(M,\xi )$ of $(W,\omega )$ corresponding to
contact $(-1)$--surgeries that undo the contact $(+1)$--surgeries
along $\LL^+$. The result will be a symplectic manifold $(W',\omega ')$
strongly filling a contact manifold $(M',\xi')$, and the latter
can be obtained from $(S^3,\xi_{st})=\partial (D^4,\omega_{st})$
by performing contact $(-1)$--surgeries (along~$\LL^-$) only.

A handlebody obtained from $(D^4,\omega_{st})$ by attaching symplectic
handles in this way is in fact a Stein filling of its
boundary contact manifold, and for those a symplectic cap had been
found earlier by Akbulut-\"Ozba\u{g}c\i\ and Lisca-Mati\'c.
The cap that fits on the Stein filling also fits on the strong
filling $(W',\omega')$, since strongly convex and strongly concave
fillings of a given contact manifold can always be glued together,
using the Liouville flow to define collar neighbourhoods of
the boundary.

\vspace{2mm}

(ii) Reduce the problem to the consideration of homology spheres only:
Let $(W,\omega )$ be a weak filling of $(M,\xi )$. We want to attach
a (weak) symplectic cobordism from $(M,\xi )$ to some integral homology sphere
$\Sigma^3$ with contact structure~$\xi'$, so as to get a weak filling
of $(\Sigma^3,\xi ')$ containing $(M,\xi )$ as a separating hypersurface.

We start from a contact surgery presentation of $(M,\xi )$ as in~(i).
For each component $L_i$ of $\LL$ we choose a Legendrian
knot $K_i$ in $(S^3,\xi_{st})$ only linked
with that component, with linking number~$1$. These $K_i$ can be chosen
in such a way that surgery with coefficient $-1$ relative to the contact framing
is the same as surgery with coefficient $0$ relative to the surface
framing. (In case you know the term: The Thurston-Bennequin invariant of
$K_i$ can be chosen to be equal to~$1$). Performing these surgeries
has the effect of killing the first integral homology.

Since $\omega$ is exact in the neighbourhood $S^1\times D^2\times
(-\varepsilon ,0]$ of a Legendrian knot in the boundary $(M,\xi )$
of $(W,\omega )$, these surgeries can be performed by attaching
symplectic $2$--handles as in the case of a strong filling. The
collection of these handles gives the desired (weak) symplectic cobordism.

\vspace{2mm}

(iii) Pass from a weak filling of a homology sphere to a strong filling:
We begin with a weak filling $(W,\omega )$
of an integral homology sphere $(\Sigma^3,\xi)$, for instance the one
obtained in~(ii); beware that we retain the original notation
for the filling. We want to
modify $\omega$ in a collar neighbourhood $\Sigma^3\times [0,1]$
of the boundary $\Sigma^3\equiv\Sigma^3\times\{ 1\}$
such that the resulting symplectic manifold is a strong filling of
the new induced contact structure on the boundary.
By (i) this can then be capped off.

Since $H^2(\Sigma^3)=0$, we can write $\omega =d\eta$ with
some $1$--form $\eta$ in a collar neighbourhood as described.
(We see that it would be enough to
have $\Sigma^3$ a rational homology sphere.) Choose a $1$--form
$\alpha$ on $\Sigma^3$ with $\xi =\ker\alpha$ and
$\alpha\wedge\omega|_{T\Sigma^3}>0$, which is possible for a weak filling.
Then set
\[ \tilde{\omega}=d(f\eta )+d(g\alpha )\]
on $\Sigma^3\times [0,1]$, where
the smooth functions $f(t)$ and $g(t)$,
$t\in [0,1]$, are
chosen as follows: Fix a small $\varepsilon >0$. Choose $f\co
[0,1]\rightarrow [0,1]$ identically $1$ on $[0,\varepsilon ]$ and
identically $0$ near~$t=1$. Choose $g\co [0,1]\rightarrow\R^+_0$
identically $0$ near $t=0$ and with $g'(t)>0$ for $t>\varepsilon /2$.

We compute
\[\tilde{\omega}=f'\, dt\wedge\eta +f\omega +g'\, dt\wedge\alpha
+g\, d\alpha,\]
whence
\begin{eqnarray*}
\tilde{\omega}^2 & = & 2ff'\, dt\wedge\eta\wedge\omega
                       + 2f'g\, dt\wedge\eta\wedge d\alpha
                       + f^2\omega^2\\
                 &   & \mbox{}
                       + 2fg'\,\omega\wedge dt\wedge\alpha
                       + 2fg\, \omega\wedge d\alpha
                       + 2gg'\, dt\wedge\alpha\wedge d\alpha.
\end{eqnarray*}
The terms appearing with the factors $f^2$, $fg'$ and $gg'$ are
positive volume forms. By choosing $g$ small on $[0,\varepsilon ]$
and $g'$ large compared with $|f'|$ and $g$
on $[\varepsilon ,1]$, one can ensure that these
positive terms dominate the three terms we cannot control.
Then $\tilde{\omega}$ is a symplectic form on the collar
and, in terms of the coordinate $s:=\log g(t)$, the symplectic
form looks like $d(e^s\alpha )$ near the boundary, with Liouville
vector field~$\partial_s$.

\section{Proof of Property P for nontrivial knots}
Here is a very rough sketch of the proof by Kronheimer and Mrowka.
It relies heavily on pretty much everything known under the sun
about gauge theory.

Let $K$ be a nontrivial knot.
It had been proved earlier by Culler-Gordon-Luecke-Shalen that
$\pi_1(K_{1/q})$ is nontrivial for
$q\not\in\{ 0,\pm 1\}$. It therefore suffices to find
a nontrivial homomorphism $\pi_1(K_1)\rightarrow\mbox{\rm SO}(3)$.

Arguing by contradiction, we assume that no such homomorphism exists.
This implies the vanishing of the instanton Floer homology group
$HF(K_1)$. By the Floer exact triangle one finds that the group
$HF(K_0)$ vanishes likewise, and so does the Fukaya-Floer homology group.

For $K$ nontrivial, results of Gabai say that $K_0$ is different from
$S^1\times S^2$ and admits a taut $2$--dimensional foliation. Eliashberg
and Thurston, in their theory of confoliations, deduce from this
the existence of a symplectic structure on $K_0\times [-1,1]$
weakly filling contact structures on the boundary components.
According to Theorem~\ref{thm:elia-filling}, by capping off these
boundaries we find a closed symplectic $4$--manifold $V$ containing $K_0$
as a separating hypersurface (and satisfying some mild cohomological
conditions).

Now, on the one hand, the Donaldson invariants of $V$ can be expressed
as a pairing on the Fukaya-Floer homology group of $K_0$
and therefore have to vanish.

On the other hand, results of Taubes say that the Seiberg-Witten invariants
of $V$ are nontrivial. By a conjecture of Witten, proved in
the relevant case by Feehan-Leness, the Donaldson invariants are likewise
nontrivial. This contradiction proves Theorem~\ref{thm:PP}.

\end{document}